\documentclass[reqno,12pt]{amsart}
\usepackage{epsfig}\usepackage{amsfonts,amsmath,amssymb,mathrsfs,verbatim,amsthm,amscd}
\numberwithin{equation}{section}
\usepackage{bm} % bold под наклоном
\newcommand{\beq}{\begin{equation}}
\newcommand{\ee}{\end{equation}}
\newcommand{\bea}{\begin{eqnarray}}
\newcommand{\eea}{\end{eqnarray}}

\usepackage{hyperref}
\def\stackreb#1#2{\ \mathrel{\mathop{#1}\limits_{#2}}}

\newcommand\id{1\kern-0.25em\text{l}}
\newcommand\zero{0\kern-0.4em\text{0}}

\newcommand{\CC}{\mathbb C}
\newcommand{\R}{\mathbb R}
\newcommand{\Z}{\mathbb Z}

\newcommand{\be}{\begin{equation}}
\newcommand{\ba}{\begin{eqnarray}}
\newcommand{\ea}{\end{eqnarray}}

\renewcommand{\epsilon}{\varepsilon}

\DeclareMathOperator{\sh}{sh}
\DeclareMathOperator{\ch}{ch}

\renewcommand{\imath}{\mathrm{i}}

\let\Re\relax

\DeclareMathOperator{\Re}{Re}

\begin{document}

\title[Complex binomial theorem and pentagon identities]
{Complex binomial theorem and pentagon identities}%

\author{N. \,M. Belousov}%

\address{N. B.:
National Research University Higher School of Economics, Moscow, Russia;
Current address: Beijing Institute of Mathematical Sciences and Applications,
Huairou district, Beijing, 101408, China}

\author{G. \,A. Sarkissian}%

\address{G. S.: Laboratory of Theoretical Physics,
JINR, Dubna, Moscow region, 141980 Russia and
Yerevan Physics Institute, Alikhanian Br. 2, 0036\, Yerevan, Armenia}

\author{V. \,P. Spiridonov}%

\address{V. S.: Laboratory of Theoretical Physics,
JINR, Dubna, Moscow region, 141980 Russia and
National Research University Higher School of Economics, Moscow, Russia}

%\vskip 1em

\begin{abstract} \noindent
We consider different pentagon identities realized by the hyperbolic hypergeometric
functions and investigate their degenerations to the level of complex hypergeometric functions.
In particular, we show that one of the degenerations yields the complex binomial theorem
which coincides with the Fourier transformation of the complex Euler beta integral {evaluation.}
At the bottom we obtain a Fourier transformation formula
for the complex gamma function. This is done with the help of a new type of the
limit $\omega_1+\omega_2\to 0$ (or $b\to \textup{i}$ in two-dimensional conformal
field theory) applied to the hyperbolic hypergeometric integrals.
\end{abstract}

\maketitle

\tableofcontents

\section{Introduction}

Investigation of hypergeometric functions over the field of complex numbers was started
by Gelfand, Graev and Vilenkin long ago \cite{GGV}. It was preceded by a seminal work of
Naimark \cite{Naimark} on the $3j$-symbols for principal series representations of the
group $SL(2,\CC)$ {(or Lorentz group)}, where the integrals over complex plane played an important role.
However, this direction of research did not get attention of
experts in the theory of special functions. Substantially, the interest to these
functions and their further investigation were inspired by the two-dimensional ($2d$)
conformal field theory (CFT)
in the middle of 1980s. Only relatively recently the group-theoretical analysis
of these functions was revived by the construction of $6j$-symbols for the group $SL(2,\CC)$
by Ismagilov \cite{Ismag2}, which was verified by a different method
in \cite{DS2017}. A description of the general structure of this class of special functions
can be found in \cite{Neretin2019}. Moreover, a proper place of these functions
was found inside the general hierarchy of special functions of hypergeometric type.
Namely, as shown in \cite{Derkachov:2021thp,Sarkissian:2020ipg},
complex hypergeometric functions emerge as special limiting forms of the
elliptic hypergeometric integrals --- the universal class of hypergeometric
special functions \cite{spi:essays}.

Degeneration of the hyperbolic hypergeometric integrals to complex
hypergeometric functions was rigorously established in \cite{Sarkissian:2020ipg}.
However, the corresponding limiting procedure was not applicable {\em per se}
in a straightforward way to hyperbolic integrals {involving exponential functions
in the integrands.} In the present paper
we continue this research line and describe further reduction of the
integrals in  \cite{Sarkissian:2020ipg} down to the simplest complex hypergeometric
forms including the complex analogue of the Fourier transformation of the Euler
gamma function.

{The most general known univariate } complex beta integral has been evaluated in \cite{Sarkissian:2020ipg}
as a limit of the elliptic beta integral and in \cite{DM2019} by a straightforward
computation. It has an interpretation as a star-triangle relation valuable for
constructing quantum integrable models. Our starting point is a
simpler complex beta integral which can be obtained as a limit from the
Faddeev-Volkov star-triangle relation \cite{FV}. Following the work
\cite{KLV}, we rewrite the corresponding identity in the form of a pentagon identity
which is used in investigations of the topology of three-dimensional manifolds.
It is well known that the hyperbolic integrals are also deeply related to $2d$ CFT
\cite{Faddeev:2000if,Ponsot:2000mt} and $3d$ supersymmetric field theories, see \cite{KLV}
and references therein.
The central charge of $2d$ quantum Liouville theory is parametrized as $c=1+6(b+b^{-1})^2$.
Therefore the limit $b\to \textup{i}$, which leads to the complex hypergeometric
functions, corresponds to $c\to 1$.

{An important} application of the hyperbolic hypergeometric integrals was found in topological
field theories
and the knot theory where they describe the so-called state integrals, or partition functions \cite{GK,KLV}.
From this point of view pentagon relations describe an invariance of the state integrals
under the $3-2$ Pachner moves. More precisely, {they correspond} to the composition
of partition functions of three ideal tetrahedra with one common edge
along which these partition functions are integrated out. This results in the partition function
of two tetrahedra glued by one common face, which corresponds to a simple product of the corresponding
partition functions (state integrals).

The main building block of the hyperbolic integrals is the Faddeev's modular
dilogarithm \cite{Fad94} {or the hyperbolic gamma function   \cite{Ruij}
(its reciprocal is known also as the Shintani's double sine function).}
We use the following notation for this function
\be
\gamma^{(2)}(u;\mathbf{\omega})= \gamma^{(2)}(u;\omega_1,\omega_2):=e^{-\frac{\pi\textup{i}}{2}
B_{2,2}(u;\mathbf{\omega}) } \gamma(u;\mathbf{\omega}),
\label{HGF}\ee
where $B_{2,2}$ is the second order multiple Bernoulli polynomial
$$
B_{2,2}(u;\mathbf{\omega})=\frac{1}{\omega_1\omega_2}
\left((u-\frac{\omega_1+\omega_2}{2})^2-\frac{\omega_1^2+\omega_2^2}{12}\right)
$$
and
\be
\gamma(u;\mathbf{\omega}):= \frac{(\tilde q e^{2\pi \textup{i} \frac{u}{\omega_1}};\tilde q)_\infty}
{(e^{2\pi \textup{i} \frac{u}{\omega_2}};q)_\infty}
=\exp\left(-\int_{\mathbb{R}+\textup{i}0}\frac{e^{ux}}
{(1-e^{\omega_1 x})(1-e^{\omega_2 x})}\frac{dx}{x}\right).
\label{int_rep}\ee
where
$$
q=e^{2\pi\textup{i}\frac{\omega_1}{\omega_2}}, \qquad \tilde {q}
=e^{-2\pi\textup{i}\frac{\omega_2}{\omega_1}},
$$ and $(a; q)_{\infty}=\prod_{k=0}^{\infty}(1-aq^k)$.

Denoting $\tau=\omega_1/\omega_2$, we see that the infinite product representation is valid for the regime $\text{Im}(\tau)>0,$ or $|q|<1$, but the integral representation defines this generalized gamma function
for $\tau\in\R$, or $|q|=1$ as well.
Since $\gamma^{(2)}(\lambda u;\lambda \omega_1,\lambda\omega_2)=\gamma^{(2)}(u;\omega_1,\omega_2)$
for any $\lambda \neq 0$, below we normalize $\omega_1\omega_2>0$ and assume that
$0\leq \arg(\omega_1)\leq \pi/2$.

This function satisfies the first order difference equations
\be\label{hp1}
{\gamma^{(2)}(u+\omega_1;\mathbf{\omega})\over \gamma^{(2)}(u;\mathbf{\omega})}=2\sin{\pi u\over \omega_2},\qquad
{\gamma^{(2)}(u+\omega_2;\mathbf{\omega})\over \gamma^{(2)}(u;\mathbf{\omega})}=2\sin{\pi u\over \omega_1},
\ee
the inversion relation $\gamma^{(2)}(u;\mathbf{\omega})\gamma^{(2)}(\omega_1+\omega_2-u;\mathbf{\omega})=1$,
and has the following asymptotics
\be\label{gamasym}
\stackreb{\lim}{u\to \infty}e^{{\pi\textup{i}\over 2}B_{2,2}(u,\omega_1,\omega_2)}\gamma^{(2)}(u;\omega_1,\omega_2)=1,
\quad {\rm for}\; {\rm arg}\;\omega_1<{\rm arg}\; u<{\rm arg}\;\omega_2+\pi,
\ee
\be\label{gamasym2}
\stackreb{\lim}{u\to \infty}e^{-{\pi\textup{i}\over 2}B_{2,2}(u,\omega_1,\omega_2)}\gamma^{(2)}(u;\omega_1,\omega_2)=1,
\quad {\rm for}\; {\rm arg}\;\omega_1-\pi<{\rm arg}\; u<{\rm arg}\;\omega_2.
\ee
$\gamma^{(2)}(u;\mathbf{\omega})$ is a meromorphic function of $u$ with the poles located
at the points
$$
u_p\in \{-n\omega_1 -m\omega_2\},\quad n,m \in \Z_{\geq0}.
$$

A general formalism for treating analogues of the quantum dilogarithm function and the corresponding
pentagon identities was described in a recent paper \cite{GK}. Namely, it was suggested to associate
such object to each local field of numbers and the case of the field of real numbers was considered in detail.
We conjecture that our consideration of the complex hypergeometric identities below
corresponds to the application of this formalism to the field of complex numbers.
Although we do not check all structural relations suggested in \cite{GK}, our key
complex hypergeometric pentagon relation
and its degenerations to simpler relations described below support this conjecture.
In particular, the bottom line of {two-dimensional} Fourier transformation for the complex gamma function
matches with the considerations of  \cite{GK}. We would like to mention also another example
of a similar general pentagon relation described in \cite{GahRos} whose proper interpretation
along the lines of  \cite{GK} is not worked out {yet}.

\section{A complex pentagon identity}

The Faddeev-Volkov star-triangle relation \cite{BMS,FV} can be written in the form
\beq\label{hyper2}
\int_{-\textup{i}\infty}^{\textup{i}\infty}\prod_{j=1}^3\gamma^{(2)}(f_j+z;\mathbf{\omega})
\gamma^{(2)}(h_j-z;\mathbf{\omega}){dz\over \textup{i}\sqrt{\omega_1\omega_2}}
=\prod_{j, k =1}^3 \gamma^{(2)}(f_j+h_k;\mathbf{\omega}),
\ee
where parameters satisfy the balancing condition $\sum_{j=1}^3 (f_j+h_j)=\omega_1+\omega_2$
and the integral converges for $\text{Re}(\omega_1+\omega_2)>0$.
The contour of integration in \eqref{hyper2} should separate two sets of poles of
the integrand  going to infinity in different directions
$$
z_{poles}\in \{h_k+n\omega_1+ m\omega_2\}\cup \{-f_k-n\omega_1- m\omega_2,\},\quad n,m \in \Z_{\geq0}.
\quad k=1,\ldots, 6.
$$
For $\text{Re}(h_k ),\text{Re}(f_k )> 0$ the contour may be chosen as the imaginary axis.
It is worth mentioning that this identity is a hyperbolic analogue of the Barnes second
lemma for a particular ordinary hypergeometric integral, whose integrand is
composed out of 6 Euler gamma functions and the result of the computation yields
a combination of 9 gamma functions.

Usually the hyperbolic hypergeometric integrals are considered under the condition
$\text{Re}(\omega_1), \text{Re}(\omega_2)>0$ (even without the assumption $\omega_1\omega_2>0$).
Their reduction to the complex hypergeometric functions takes place in the singular
limit $\omega_1+\omega_2\to 0$. For taking such a limit we use the following parametrization
fitting our conventions
\beq\label{omi}
\omega_1=\textup{i}+\delta, \qquad \omega_2=-\textup{i}+\delta,\quad \delta\to 0^+,
\ee
which corresponds to the choice
\beq\label{om1om22}
b=\sqrt{\omega_1\over \omega_2}=\textup{i}+\delta+O(\delta^2).
\ee
As a result, the function $\gamma^{(2)}(u;\mathbf{\omega})$ starts to blow up at special
points of the argument $u$ lying on the imaginary axis. The exact asymptotics was
rigorously worked out in \cite{Sarkissian:2020ipg}--- one has uniformly on the compacta
\beq\label{gam2lim2}
\gamma^{(2)}(\textup{i}\sqrt{\omega_1\omega_2}(n+x\delta);\omega_1,\omega_2)\stackreb{=}{\delta\to 0^+} e^{\frac{\pi \textup{i}}{2}n^2} (4\pi\delta)^{\textup{i}x-1}{\bf \Gamma}(x,n),
\ee
where $n\in \Z, \, x\in\CC$. Heuristically, this limit was considered earlier in \cite{BMS}.
Here ${\bf \Gamma}(x,n)$ is the gamma function over the field of complex numbers
\begin{equation}
{\bf \Gamma}(x,n)\equiv{\bf \Gamma}(\alpha|\alpha'):=\frac{\Gamma(\alpha)}{\Gamma(1-\alpha')}
=\frac{\Gamma(\frac{n+\textup{i}x}{2})}{\Gamma(1+\frac{n-\textup{i}x}{2})},
\quad \alpha=\frac{n+\textup{i}x}{2},\; \alpha'=\frac{-n+\textup{i} x}{2},
\label{Cgamma}\end{equation}
where $x\in \CC$ and $n\in\Z$.
The well known relation for the Euler gamma function $\Gamma(x)\Gamma(1-x)=\pi/\sin\pi x$
leads to the identities
\beq
{\bf \Gamma}(\alpha|\alpha') =(-1)^{\alpha-\alpha'}{\bf \Gamma}(\alpha'|\alpha), \qquad
{\bf \Gamma}(x,-n)=(-1)^n{\bf \Gamma}(x,n),
\label{reflCgamma0}\ee
and
\beq
{\bf \Gamma}(\alpha|\alpha'){\bf \Gamma}(1-\alpha|1-\alpha')  =(-1)^{\alpha-\alpha'}, \qquad
{\bf \Gamma}(x,n){\bf \Gamma}(-x-2\textup{i},n)=1.
\label{reflCgamma}\ee
Functional equations have the form
\begin{eqnarray*} &&
{\bf \Gamma}(\alpha+1|\alpha') ={\bf \Gamma}(x-\textup{i},n+1)=\alpha{\bf \Gamma}(\alpha|\alpha'),\quad
\\ &&
{\bf \Gamma}(\alpha|\alpha'+1) ={\bf \Gamma}(x-\textup{i},n-1)=-\alpha' {\bf \Gamma}(\alpha|\alpha').
\end{eqnarray*}

Now one sets in (\ref{hyper2})
\beq
f_j=\textup{i}\sqrt{\omega_1\omega_2}(N_j+\delta a_j),
\qquad h_j=\textup{i}\sqrt{\omega_1\omega_2}(M_j+\delta b_j),\quad j=1,2,3,
\ee
and the integration variable is replaced by an analogous expression
$z=\textup{i}\sqrt{\omega_1\omega_2}(n+\delta y)$, where $n, N_j, M_j \in \Z.$
As a result, in the limit $\delta\to 0$ the balancing condition
splits into two relations
$$
\sum_{j=1}^3( N_j+M_j)=0,\qquad \sum_{j=1}^3( a_j+b_j)=-2\textup{i} .
$$
According to \cite{Sarkissian:2020ipg}, in the limit $\delta\to 0^+$ relation (\ref{hyper2}) reduces to
the identity
\bea \nonumber &&
\frac{1}{4\pi} \sum_{n\in \Z}
\int_{-\infty}^{\infty}(-1)^{n+\sum_{j=1}^3M_j}\prod_{j=1}^3{\bf \Gamma}(a_j+y,N_j +n)
{\bf \Gamma}(b_j-y,M_j -n)dy
\\ && \makebox[6em]{}
 =\prod_{j,k=1}^3{\bf \Gamma}(a_j+b_k,N_j+M_k).
 \label{limtrafoI}\eea
In \cite{DMV2017} it was shown that this formula represents the Mellin-Barnes form of the star-triangle relation
directly related to the complex beta integral of \cite{GGV}.
A brute force method of computing this integral was presented also in \cite{DM2019}.
{We would like to mention also} that this identity is a complex analogue of the second Barnes lemma.

Let us denote
$$
B(x,y) :=\frac{\gamma^{(2)}(x,y;\mathbf{\omega})}{\gamma^{(2)}(x + y;\mathbf{\omega})}.
$$
Then, as shown in \cite{KLV}, the exact integration formula \eqref{hyper2} can be rewritten
in the form of pentagonal (five term) relation
\beq\label{pent_hyp}
\int_{-\textup{i}\infty}^{\textup{i}\infty}\prod_{j=1}^3
B(f_j+z,h_j-z){dz\over \textup{i}\sqrt{\omega_1\omega_2}}
= B(f_2 + h_1, f_3 + h_2)B(f_1 + h_2,f_3 + h_1).
\ee

Our main starting point is the following complex pentagon identity
\begin{eqnarray} \label{pent_comp}  &&
\sum_{n\in\Z}\int_{-\infty}^{\infty}(-1)^{n+\sum_{j=1}^3M_j}\prod_{j=1}^3
B(a_j+y,N_j+n; b_j-y,M_j-n)\frac{dy}{4\pi}
\\ && \makebox[1em]{}
= B(a_2 + b_1,N_2+M_1; a_3 + b_2,N_3+M_2)B(a_1 + b_2,N_1+M_2;a_3 +b_1,N_3+M_1),
\nonumber\end{eqnarray}
where
$$
B(x,n;y,m):=\frac{{\bf \Gamma}(x,n){\bf \Gamma}(y,m)}{{\bf \Gamma}(x+y,-n-m)}.
$$
It is obtained from relation \eqref{limtrafoI} in the same way as identity \eqref{pent_hyp}
follows from the hyperbolic beta integral \eqref{hyper2}.
Notice that in \eqref{pent_comp} not all external parameters $N_j$ (or $M_j$) are hidden inside the composite
$B$-functions. This means that there must be some specific interpretation of the discrete variables $N_j,\, M_j$,
since the cases of even and odd integer numbers $\sum_{j=1}^3N_j$ yield different pentagon identities.

Consider degeneration limits of the hyperbolic pentagon identity \eqref{pent_hyp}.
Such degenerations of general hyperbolic integrals were considered in detail
in \cite{bult,rai:limits}. Here we follow
the paper \cite{Sarkissian:2020ipg} and extend its analysis by considering integrals involving
exponential functions in the integrands.
The limit  $f_3\to +\textup{i}\infty, \, h_3\to -\textup{i}\infty$ with $f_3+h_3$ fixed, simplifies
the {exact integration formula \eqref{hyper2}} to
\begin{eqnarray}\nonumber  &&
\int_{-\textup{i}\infty}^{\textup{i}\infty}
e^{\frac{\pi \textup{i} z}{\omega_1\omega_2}(f_1+f_2+h_1+h_2)}
\prod _{j=1}^2\gamma^{(2)}(f_j+z;\mathbf{\omega})
\gamma^{(2)}(h_j-z;\mathbf{\omega}){dz\over \textup{i}\sqrt{\omega_1\omega_2}}
\\ && \makebox[2em]{}
=e^{{\pi \textup{i}\over \omega_1\omega_2}(h_1h_2 - f_1f_2)}
\frac{\prod_{j,k=1}^2\gamma^{(2)}(f_j+h_k;\mathbf{\omega})}
{\gamma^{(2)}(\sum_{j=1}^2(f_j+h_j);\mathbf{\omega})}.
\label{hyp4-5}\end{eqnarray}
The integral on the left-hand side converges in the upper limit for $\text{Re}(\omega_1+\omega_2)>0$
and in the lower limit for $\text{Re}(\omega_1+\omega_2)>\text{Re}(f_1+f_2+h_1+h_2)$.
In  \cite{SaSpStek} {such a} relation was obtained in a much more general context
as a hyperbolic beta integral on the general lens space and an extension of \eqref{hyp4-5}
to the multivariate functions on the root systems was described in \cite{DM2019,SaSpJGP}.
{The derived identity} is a hyperbolic analogue of the first Barnes lemma, which converts an integral
of a product of 4 Euler gamma functions into a combination of 5 gamma functions. Note that it
can be rewritten in the following form
\begin{eqnarray}\nonumber &&
\int_{-\textup{i}\infty}^{\textup{i}\infty}
e^{\frac{\pi \textup{i} z}{\omega_1\omega_2}(f_1+f_2+h_1+h_2)}
\prod _{j=1}^2B(f_j+z,h_j-z) {dz\over \textup{i}\sqrt{\omega_1\omega_2}}
\\ && \makebox[2em]{}
=e^{{\pi \textup{i}\over \omega_1\omega_2}(h_1h_2 - f_1f_2)}B(f_1+h_2,f_2+h_1),
\label{hyp4-5'}\end{eqnarray}
which resembles a relation between partition functions for a 2--1 move of some geometrical figures.

Applying the limit $b\to \textup{i}$ in a straightforward way as described above, we see that
the exponentials in the left- and right-hand side expressions disappear and there emerges
the equality \cite{DM2019,SaSpJGP}
\begin{eqnarray} \nonumber &&
\frac{1}{4\pi}\sum_{m\in \Z}\int_{\R }
 \prod_{\ell=1}^{2}{\bf \Gamma}(a_\ell+ x,N_\ell+ m){\bf \Gamma}(b_\ell- x,M_\ell- m)dx
\\  && \makebox[2em]{}
=
{\prod_{\ell,s=1}^{2}{\bf \Gamma}(a_\ell+b_s,N_\ell+M_s)\over {\bf \Gamma}(\sum_{s=1}^{2}(a_s+b_s),\sum_{s=1}^{2}(N_s+M_s))},
\label{complexAI2}\end{eqnarray}
or
\begin{eqnarray}\nonumber &&
\sum_{m\in \Z}\int_{-\infty}^{\infty}
\prod _{j=1}^2B(a_j+x,N_j+m;b_j-x,M_j-m)\frac{dx}{4\pi}
\\ && \makebox[2em]{}
=B(a_1+b_2,N_1+M_2;a_2+b_1,N_2+M_1).
\label{comp4-5'}\end{eqnarray}

Another option for reducing identity \eqref{hyper2} to a lower level consists in taking
simultaneously two parameters, say, $f_2, f_3\to +\textup{i}\infty$ and one parameter $h_3\to -\textup{i}\infty$
in such a way that the sum $f_2+f_3+h_3$ remains a finite fixed quantity. This leads to the
identity corresponding to the choice of `+' sign in the exponentials of the following relation
\begin{eqnarray}\nonumber &&
\int_{-\textup{i}\infty}^{\textup{i}\infty}
e^{\pm \frac{\pi \textup{i} z}{2\omega_1\omega_2}( \omega_1+\omega_2 +2(f_1+h_1+h_2)-z)}
\gamma^{(2)}(f_1+z,h_1-z,h_2-z;\mathbf{\omega}){dz\over \textup{i}\sqrt{\omega_1\omega_2}}
\\ && \makebox[2em]{}
=e^{\pm {\pi \textup{i}\over 2\omega_1\omega_2}(f_1^2+2f_1(h_1+h_2-\frac{\omega_1+\omega_2}{2}) +4h_1h_2)}
\gamma^{(2)}(f_1+h_1,f_1+h_2;\mathbf{\omega}) .
\label{hyp3-2}\end{eqnarray}
The same identity arises in the limit $f_2\to +\textup{i}\infty$ applied to relation \eqref{hyp4-5}.
After taking the limit $h_2\to -\textup{i}\infty$ in \eqref{hyp4-5}, reflecting the integration
variable $z\to -z$ and changing notation $(f_1,f_2,h_1)\to(h_1,h_2,f_1)$ one obtains relation
\eqref{hyp3-2} with `$-$' sign in the exponentials.
For `+' sign, the integral on the left-hand side converges in the upper limit for $\text{Re}(\omega_1+\omega_2)>0$
and in the lower limit for $\text{Re}(\omega_1+\omega_2)>2\text{Re}(f_1+h_1+h_2)$. For `$-$' sign,
the convergence conditions in the upper and lower limits flip.

Equality \eqref{hyp4-5} can be simplified further by taking the limit
 $f_2\to -\textup{i}\infty$ and $h_2\to \textup{i}\infty$ with $f_2+h_2=\alpha$ kept fixed, which yields
\bea\label{namer2} &&
\int_{-\textup{i}\infty}^{\textup{i}\infty}e^{{\pi \textup{i}z\over \omega_1\omega_2 }
(2\alpha+ h+f-\omega_1-\omega_2)}
\gamma^{(2)}(f+z,h-z;\mathbf{\omega}) \frac{dz}{\textup{i}\sqrt{\omega_1\omega_2}}
\\ \nonumber  &&  \makebox[4em]{}
=e^{{\pi \textup{i}\over 2\omega_1\omega_2 }
(h-f)(2\alpha+h+f-\omega_1-\omega_2)} \gamma^{(2)}(\alpha,\omega_1+\omega_2-h-f-\alpha,f+h;\mathbf{\omega}),
\eea
where we denoted $f_1=f$ and $h_1=h$.
The integral on the left-hand side converges for in the upper limit for $\text{Re}(\alpha)>0$
and in the lower limit for $\text{Re}(\alpha)<\text{Re}(\omega_1+\omega_2 -h-f)$.
This is again a kind of a pentagon identity, when two hyperbolic gamma functions are integrated
to the product of three such functions. However, we can rewrite this identity in an alternative
a kind of 1--1 move form. Denote
$$
\lambda:= \alpha+\tfrac12(h+f-\omega_1-\omega_2),\quad g:=\omega_1+\omega_2-h-f.
$$
Then, shifting the integration variable $z\to z+\tfrac12(h-f)$, we deform appropriately
the integration contour (for $\text{Re}(h),\text{Re}(f)> 0$ this is just a shift of the imaginary axis)
and write
\bea\label{pentafkv5}  \makebox[-2em]{}
\int_{-\textup{i}\infty}^{\textup{i}\infty}
e^{2\pi\textup{i}\lambda z\over \omega_1\omega_2}
 \gamma^{(2)}\left(\pm z+{\bar{g}\over 2};\mathbf{\omega}\right)
{dz\over \textup{i}\sqrt{\omega_1\omega_2}}=\gamma^{(2)}(\bar{g};\mathbf{\omega})
\gamma^{(2)}\left(\pm \lambda+{g\over 2};\mathbf{\omega}\right),
\eea
where $\bar g=\omega_1+\omega_2-g$ and we use the notation
$$ 	f(x \pm y) = f(x + y) f(x - y). $$
This integral converges under the condition that $\text{Re}(g)>2|\text{Re}(\lambda)|\geq 0$.
Equivalently, {for purely imaginary $\lambda$  it can be written as
a Fourier transformation
\begin{eqnarray}&& \makebox[-2em]{}
\int_{-\textup{i}\infty}^{\textup{i}\infty}
e^{{2\pi \textup{i}\lambda z\over \omega_1\omega_2 } }
\tilde B(\tfrac{\bar g}{2}+z,\tfrac{\bar g}{2}-z) {dz\over \textup{i}\sqrt{\omega_1\omega_2}}
=\tilde B(\tfrac{g}{2}+\lambda,\tfrac{g}{2}-\lambda)
\label{hyp1-1}\end{eqnarray}
of {the } Boltzmann weight}
$$
\tilde B(x,y)=\frac{\gamma^{(2)}(x;\mathbf{\omega})\gamma^{(2)}(y;\mathbf{\omega})}
{\sqrt{\gamma^{(2)}(x+y;\mathbf{\omega})}}.
$$

Relation \eqref{pentafkv5} is obtained as a degeneration of the general pentagon identity \eqref{pent_hyp}
and represents again a 5-term relation describing a kind of  2--3 move. Therefore we call it the
reduced pentagon identity. It was explicitly written in \cite{Faddeev:2000if,Ponsot:2000mt},
and, in particular, in \cite{Faddeev:2000if} it was used in order to prove the
following operator pentagon identity
\begin{equation}
e_b(P)e_b(X)=e_b(X)e_b(P + X)e_b(P), \qquad [P,X]=\frac{1}{2\pi\textup{i}},
\label{oper_pent}\end{equation}
where
$$
e_b(z)=\frac{1}{\gamma(\frac{\omega_1+\omega_2}{2}-\textup{i}\sqrt{\omega_1\omega_2}z;\mathbf{\omega})}.
$$
Similar operator identities for more general functional relations up to \eqref{hyper2}
were discussed in \cite{Kashaev,Volkov}.

Equation \eqref{namer2} can be simplified further in two different ways. First, we take the limit
$h\to -\textup{i}\infty$. Then we shift the integration variable $z\to z-f$ and deform appropriately the
integration contour, which removes the parameter $f$ completely from the integral. This yields
the following exact formula
\bea\label{namer3} &&
\int_{\textup{i}\mathbb{R}+0}
e^{\frac{2\pi \textup{i}\lambda z}{\omega_1\omega_2} + \frac{\pi \textup{i}}{2}(B_{2,2}(z;\mathbf{\omega})
- B_{2,2}(0;\mathbf{\omega}))}
\gamma^{(2)}(z;\mathbf{\omega})\frac{dz}{\textup{i}\sqrt{\omega_1\omega_2}}
=e^{-{\pi \textup{i}\over 2} B_{2,2}(\lambda;\mathbf{\omega})}
\gamma^{(2)}(\lambda;\mathbf{\omega}),
\eea
where we replaced the symbol $\alpha$ by $\lambda$ in order to make the formula
close in notation to \eqref{hyp1-1}.
It looks like the Fourier transformation of the function \eqref{int_rep}, however
the integral converges under the condition $0<\text{Re}(\lambda)< \tfrac12 \text{Re}(\omega_1+\omega_2)$,
i.e. the standard requirement $\text{Re}(\lambda)=0$ is not applicable.
Note also that the contour of integration should go around the singular point $x=0$ from the right.

Second, we take in the equality \eqref{namer2} the limit $f\to +\textup{i}\infty$. Then we
shift the integration variable $z\to z+h$, deform appropriately the
integration contour, and reflect the resulting integration variable $z\to -z$.
As a result, we obtain the relation
\bea\label{namer4} &&
\int_{\textup{i}\mathbb{R}+0}
e^{-\frac{2\pi \textup{i}\lambda z}{\omega_1\omega_2} - \frac{\pi \textup{i}}{2}(B_{2,2}(z;\mathbf{\omega})
- B_{2,2}(0;\mathbf{\omega}))}
\gamma^{(2)}(z;\mathbf{\omega})\frac{dz}{\textup{i}\sqrt{\omega_1\omega_2}}
=e^{{\pi \textup{i}\over 2} B_{2,2}(\lambda;\mathbf{\omega})}
\gamma^{(2)}(\lambda;\mathbf{\omega}),
\eea
where again the integral converges for $0<\text{Re}(\lambda)< \tfrac12 \text{Re}(\omega_1+\omega_2)$.
This equality looks like the inverse Fourier transformation of relation \eqref{namer3}, but
the standard proof of this inversion does not go through. We cannot multiply the left-hand side
of \eqref{namer3} by $e^{-\frac{2\pi \textup{i}\lambda y}{\omega_1\omega_2}}$ and
integrate over $\lambda$ along the imaginary axis because of the integral divergency.
Still, with a slight modification of the inverse Fourier transformation one can derive
one formula from another.

The key results of the present note consist in the derivation of complex analogues
of the relations \eqref{pentafkv5} and  \eqref{namer3}.
We shall obtain them by applying an extension of the limiting procedure $b\to \textup{i}$
of \cite{Sarkissian:2020ipg} described above.
Detailed discussion of similar limits for \eqref{hyp4-5} and \eqref{hyp3-2}
is postponed to a future work.

\section{A new degeneration limit}

To investigate the degeneration of above hyperbolic integrals to the complex hypergeometric level we need
a new limiting formula for the following ratio of two hyperbolic gamma functions
\begin{equation}\label{g-ratio'}
\frac{\gamma^{(2)}(z + y;\mathbf{\omega})} {\gamma^{(2)} (z;\mathbf{\omega})}=
\frac{e^{\frac{\pi \imath}{2}B_{2,2}(z)}(e^{2\pi\imath\frac{z}{\omega_2}};q)_\infty
(\tilde q e^{2\pi\imath\frac{z+y}{\omega_1}};\tilde q)_\infty}
{e^{\frac{\pi \imath}{2} B_{2,2}(z + y)}(\tilde q e^{2\pi\imath\frac{z}{\omega_1}};\tilde q)_\infty
(e^{2\pi\imath\frac{z+y}{\omega_2}};q)_\infty}.
\ee
Recall the standard $q$-analogue of the $_1F_0$-hypergeometric function
\begin{align}
	_1\varphi_0(a;q;v) =\sum_{n=0}^{\infty}\frac{(a; q)_n}{(q; q)_n}v^n, \quad |v|<1, \quad |q|<1,
\end{align}
where $(a; q)_n=\prod_{k=0}^{n-1}(1-aq^k)$, and the $q$-binomial theorem \cite{GR}
\begin{align}
	_1\varphi_0(a;q;v)= \frac{(av; q)_{\infty}}{(v; q)_{\infty}}.
\end{align}

Then, one can identify in \eqref{g-ratio'} ratios of infinite products with the base $q$
and $\tilde q$ with two $_1\varphi_0$-series with different bases and write
\begin{align}\label{g-ratio}
		\frac{\gamma^{(2)}(z + y;\mathbf{\omega})} {\gamma^{(2)} (z;\mathbf{\omega})}
=  \frac{e^{\frac{\pi \imath}{2} B_{2,2}(z)}\,
{}_1\varphi_0\bigl(e^{2\pi \imath \frac{y}{\omega_1}}; \tilde q;  \tilde q e^{2\pi \imath \frac{z}{\omega_1}} \bigr)}
{e^{\frac{\pi \imath}{2} B_{2,2}(z+y)}\, {}_1\varphi_0\bigl(e^{2\pi \imath \frac{y}{\omega_2}}; q; e^{2\pi \imath \frac{z}{\omega_2}} \bigr)}.
\end{align}
Below we shall use the degeneration to the standard binomial theorem
\begin{align}\label{gqt-lim}
{}_1\varphi_0(q^c;q;v)\stackreb{=}{q\to 1^-}  {}_1F_0(c;v)=(1-v)^{-c}.
\end{align}

Now we insert into \eqref{g-ratio} the parametrizations \eqref{omi} and
\begin{align}\label{zy_par}
z=\imath \sqrt{\omega_1 \omega_2} \bigl(N + \beta \bigr), \quad
y =  \imath \sqrt{\omega_1 \omega_2} (m + u \delta),\quad N, m\in\Z,\; \beta, u\in\CC,
\end{align}
and take the limits $\delta\to 0^+$, $N\to \pm \infty$ in such a way that
\begin{align}\label{alpha}
	N \delta = \alpha \in \mathbb{R}
\end{align}
is a fixed constant.
In this limit we have the following asymptotic estimates
\begin{align}\label{exp}
	\begin{aligned}
	& q=e^{2\pi \imath \frac{\omega_1}{\omega_2} } = e^{- 4\pi \delta + O(\delta^2)},
&& \quad \tilde q=e^{-2\pi \imath \frac{\omega_2}{\omega_1} } = e^{- 4\pi \delta + O(\delta^2)},
\\[6pt]
	& e^{2\pi \imath \frac{z}{\omega_2}} = e^{-2\pi (\alpha +\imath \beta) + O(\delta)},
&& \quad e^{2\pi \imath \frac{z}{\omega_1}} = e^{- 2\pi(\alpha - \imath \beta) + O(\delta)},
\\[6pt]
&  e^{2 \pi \imath \frac{y}{\omega_2}}
% = e^{- 2\pi (m + \imath u) \delta + O(\delta^2)}
=q^{\frac{m + \imath u}{2}+O(\delta)}, \qquad
&& \quad e^{2 \pi \imath \frac{y}{\omega_1}}
% = e^{- 2\pi (m - \imath u) \delta + O(\delta^2)}
=\tilde q^{\frac{m - \imath u}{2}+O(\delta)}.
	\end{aligned}
\end{align}
As we see, both $_1\varphi_0$-series in \eqref{g-ratio} simultaneously go to $_1F_0$-series
yielding the powers of elementary functions according to \eqref{gqt-lim}.
Calculation of the asymptotics for the exponentials of Bernoulli polynomials yields
\begin{equation}\label{Ber_lim_even}
	e^{\frac{\pi \imath}{2} [ B_{2,2}(z) - B_{2,2}(z + y) ]}
\stackreb{=}{\delta\to 0^+, \, |N|\to \infty} e^{\pi\imath Nm}e^{\frac{\pi \imath}{2} m^2
+ \pi\imath( \alpha u +  \beta m) + O(\delta)},
\end{equation}
i.e. the limit is not well defined for odd $m$.

Combining these results we come to the key limiting formula needed to us
\begin{eqnarray}\label{lim2} &&
e^{\pi\imath Nm}\frac{\gamma^{(2)} \bigl( \imath \sqrt{\omega_1 \omega_2} \bigl[ N + \beta + m + u \delta \bigr] \bigr)}{\gamma^{(2)} \bigl( \imath \sqrt{\omega_1 \omega_2} \bigl[ N + \beta\bigr] \bigr)}
\\ && \makebox[2em]{}
	\underset{\substack{\; \delta \to 0^+,\, |N| \to \infty \\[2pt] N\delta=\alpha } }{=} e^{\frac{\pi \imath}{2}m^2}\, \biggl( 2 \sh \pi(\alpha + \imath \beta) \biggr)^{\frac{m + \imath u}{2} } \biggl( 2 \sh \pi(\alpha - \imath \beta) \biggr)^{\frac{- m + \imath u}{2} }
\nonumber\end{eqnarray}
where $N, m \in \mathbb{Z}$, $\alpha\in\R$, $\beta, u \in \mathbb{C}$. Note that for the
complex variable $z=\sh\pi(\alpha+\imath\beta)$ we have in general branch singularities which
disappear only for $\beta\in\R$.
The fact that we have different limits for the ratio of hyperbolic
gamma functions along two different discrete sets of points
$$
\delta^{(k)}_1=\frac{\alpha}{2k}\stackreb{=}{k\to\pm\infty} 0^+, \quad
\delta^{(k)}_2=\frac{\alpha}{2k+1}\stackreb{=}{k\to\pm\infty} 0^+
$$
is an interesting artifact of our algorithm.
This delicate point should be taken into account whenever one applies such a limit to a given
relation. Above we resolved it by introducing in \eqref{lim2} an artificial compensating
sign function $(-1)^{Nm}$. Thus, we obtained a new $b\to\imath$ asymptotic formula for $\gamma^{(2)}$-functions
when their arguments behave in a particular way, which is different from the limit \eqref{gam2lim2}
rigorously established in \cite{Sarkissian:2020ipg}.

\section{Degeneration of a pentagon identity to the complex binomial theorem}

The $b\to 0$ limit of the reduced pentagon identity \eqref{pentafkv5} is considered in detail
in \cite{BDKK2}. In this section we consider its complementary limit $b\to \textup{i}$.
For that we first shift the integration variable $z=y-\bar g/2$ and rewrite it as follows
\bea\label{pentafkv5'}  \makebox[-2em]{}
\frac{1}{\gamma^{(2)}(\bar{g};\mathbf{\omega})}\int_{\textup{i}\R+0}
e^{\frac{2\pi\textup{i}\lambda y}{\omega_1\omega_2}}
 \gamma^{(2)}\left(y,\bar{g}-y;\mathbf{\omega}\right)
{\frac{dy}{\textup{i}\sqrt{\omega_1\omega_2}}}=e^{\frac{\pi\textup{i}\lambda \bar g}{\omega_1\omega_2} }
\gamma^{(2)}\left(\pm \lambda+{g\over 2};\mathbf{\omega}\right),
\eea
where the contour passes the $y=0$ point from the right. Now we parametrize
\begin{align} \label{g}
	& \bar{g} = \imath  \sqrt{\omega_1 \omega_2}(r + h \delta), && r \in \mathbb{Z}, \qquad h \in \mathbb{C}, \\[6pt]
	& y = \imath\sqrt{\omega_1 \omega_2}(m + u \delta), && m \in \mathbb{Z}, \qquad u \in \mathbb{R}, \\[6pt] \label{l}
	& \lambda = \imath\sqrt{\omega_1 \omega_2} \big(N+\beta \big), && N \in \mathbb{Z},\qquad \beta\in\mathbb{C},
\end{align}
and take the limits $\delta \to 0^+,\, N\to\pm \infty$ with the assumption that
\begin{equation}
	N \delta = \alpha\in\R.
\end{equation}

The integrand in \eqref{pentafkv5'} has rich set of poles which approach for $\delta\to 0^+$ the integration contour
(the imaginary axis) and pinch it  at infinitely many points $y=\imath  m,\, m\in\Z$.
The asymptotics of $\gamma^{(2)}(y;\mathbf{\omega})$ in this regime is given by formula \eqref{gam2lim2}.
In \cite{Sarkissian:2020ipg} it was shown that for a function $\Delta(y)$ formed from the $\gamma^{(2)}(g_k+y;\mathbf{\omega})$-functions with different parameters $g_k$ constrained
similar to $\bar g$ above, the asymptotics of its integral over the imaginary axis for $\delta\to 0^+$
has the form
\beq\label{intdel}
\int_{-\textup{i}\infty}^{\textup{i}\infty}\Delta(y){dy\over \textup{i}\sqrt{\omega_1\omega_2}}
\stackreb{=}{\delta\to 0^+}
\sum_{m\in \Z}
\int_{-\infty}^{\infty} \left[
\stackreb{\lim}{\delta\to 0}\delta\Delta(\textup{i}\sqrt{\omega_1\omega_2}(m+u\delta))\right] du,
\ee
provided on the right-hand side one has uniform asymptotic limit of the integrand yielding convergent
sums of integrals times a power of $\delta$.

In this way we obtain
$$
e^{\frac{2\pi\textup{i}\lambda y}{\omega_1\omega_2} }
=e^{-2\pi\textup{i}(N+\beta)(m+u\delta)}\stackreb{=}{\delta\to 0^+, \, N \to\infty}
e^{-2\pi\textup{i}(\alpha u+\beta m)}
$$
and
$$
\frac{\gamma^{(2)}\left(y,\bar{g}-y;\mathbf{\omega}\right)}{\gamma^{(2)}(\bar{g};\mathbf{\omega})}
\stackreb{=}{\delta\to 0^+} \frac{e^{\pi\imath m(1-r)}{\bf\Gamma}(u,m){\bf\Gamma}(h-u,r-m)}{4\pi\delta\, {\bf\Gamma}(h,r)},
$$
where we used the relation $e^{\pi\imath m^2}=e^{\pi\imath m}$.
Therefore we have
\bea
&&\lim_{\delta\to0+} {1\over\gamma^{(2)}(\bar{g};\mathbf{\omega}) } \int_{-\textup{i}\infty}^{\textup{i}\infty}
e^{\frac{2\pi\textup{i}y\lambda}{\omega_1\omega_2}}
 \gamma^{(2)}\left(y,\bar{g}-y;\mathbf{\omega}\right)
{\frac{dy}{\textup{i}\sqrt{\omega_1\omega_2}}}
 \nonumber\\
&& =  \frac{1}{4\pi \bm{\Gamma}(h, r)}
 \sum_{m \in \mathbb{Z}} \int_{\mathbb{R}-\imath 0}  e^{\pi\textup{i} m(1-r)}
 e^{-2\pi \imath(  \alpha u +  m\beta) } \; {\bf\Gamma}(u,m) {\bf\Gamma}(h-u, r- m)\, du.	
\eea

As to the right-hand side of \eqref{pentafkv5'}, we have for the exponential factor
\begin{equation}\label{prefact}
e^{\frac{\pi\textup{i}\lambda \bar g}{\omega_1\omega_2} }
=e^{-\pi\textup{i}(N+\beta)(r+h\delta)}\stackreb{=}{\delta\to 0^+, \, |N| \to\infty}
(-1)^{Nr}e^{-\pi\textup{i}(\alpha h+\beta r)},
\end{equation}
which for odd $r$ has a sign alternating multiplier $(-1)^N$, i.e. the limit is not well defined.
However, this drawback is cured by the asymptotics of the rest expression on the right-hand side.
Indeed,
\be\label{intfo}
\gamma^{(2)}\left(\pm \lambda+{g\over 2}\right)={\gamma^{(2)}\left(\lambda+{\omega_1+\omega_2-\bar{g}\over 2}\right)\over \gamma^{(2)}\left( \lambda+{\omega_1+\omega_2+\bar{g}\over 2}\right)}
={\gamma^{(2)}\left(z-\bar{g}\right)\over \gamma^{(2)}\left( z\right)},
\ee
where
\be
z=\lambda+{\omega_1+\omega_2+\bar{g}\over 2}=\imath\sqrt{\omega_1 \omega_2} \bigl( N + \beta +\tfrac12 r
+ (\tfrac12 h-\imath ) \delta\bigr)+O(\delta^2).
\ee

Applying relation \eqref{lim2} to \eqref{intfo},  we obtain
\bea && \makebox[-2em]{}
\lim_{\delta\to 0^+, \, N \to\infty}\gamma^{(2)}\left(\pm \lambda+{g\over 2}\right)
\nonumber \\ &&
=(-1)^{Nr} e^{\frac{\pi \imath r^2}{2}}\,
\bigl( 2 \sh \pi(\alpha + \imath(\beta+\tfrac12 r)) \bigr)^{\frac{-r - \imath h}{2}}
\bigl( 2 \sh \pi(\alpha - \imath (\beta+\tfrac12 r)) \bigr)^{\frac{r - \imath h}{2} }.
\eea
Combining this expression with \eqref{prefact} we see that the troubling multiplier $(-1)^{Nr}$
cancels out and we obtain the formula
\bea && \makebox[-2em]{}
\frac{1}{4\pi \bm{\Gamma}(h, r)}
 \sum_{m \in \mathbb{Z}} \int_{\mathbb{R}-\imath 0} e^{\pi\textup{i} m(1-r)}
 e^{-2\pi \imath(  \alpha u +  m\beta) } \; {\bf\Gamma}(u,m) {\bf\Gamma}(h-u, r- m)\,  du
\nonumber \\ &&
= e^{\pi\textup{i}(\frac{1}{2}r^2-\alpha h-\beta r)}\,
\bigl( 2 \sh \pi(\alpha + \imath(\beta+\tfrac12 r)) \bigr)^{\frac{-r - \imath h}{2}}
\bigl( 2 \sh \pi(\alpha - \imath (\beta+\tfrac12 r)) \bigr)^{\frac{r - \imath h}{2} }.
\label{qbin_comp}\eea

Shifting $\beta\to \beta-r/2$ we simplify this expression to
\bea &&
\frac{1}{4\pi \bm{\Gamma}(h, r)}
 \sum_{m \in \mathbb{Z}} e^{\pi\imath(m-r)} \int_{\mathbb{R}-\imath 0}
 e^{-2\pi \imath(  \alpha (u-\frac12 h) +  \beta(m-\frac12 r)) } \; {\bf\Gamma}(u,m) {\bf\Gamma}(h-u, r- m)\, du
\nonumber \\ && \makebox[4em]{}
= \bigl( 2 \sh \pi(\alpha + \imath\beta) \bigr)^{\frac{-r - \imath h}{2}}
\bigl( 2 \sh \pi(\alpha - \imath \beta) \bigr)^{\frac{r - \imath h}{2} }.
\eea

Let $r=2\ell$ is even and $h=2s$. Then we can shift $m\to m+\ell$, $u\to u+s$  and obtain
\bea &&\label{even_l}
\frac{1}{4\pi \bm{\Gamma}(2s, 2\ell)}
 \sum_{m \in \mathbb{Z}} e^{\pi\imath(m-\ell)} \int_{\mathbb{R}}
 e^{-2\pi \imath(  \alpha u +  \beta m) } \; {\bf\Gamma}(s \pm u,\ell \pm m)\, du
\nonumber \\ && \makebox[4em]{}
= \bigl( 2 \sh \pi(\alpha + \imath\beta) \bigr)^{-\ell - \imath s}
\bigl( 2 \sh \pi(\alpha - \imath \beta) \bigr)^{\ell - \imath s}.
\eea

Let $r=2\ell+1$ is odd and $h=2s$. Then we can shift again $m\to m+\ell$, $u\to u+s$  and obtain
\bea && \makebox[-2em]{}
\frac{1}{4\pi \bm{\Gamma}(2s, 2\ell+1)}
 \sum_{m \in \mathbb{Z}} e^{\pi\imath(m-\ell-1)} \int_{\mathbb{R}}
 e^{-2\pi \imath(  \alpha u +  \beta (m-\frac12)) } \; {\bf\Gamma}(u+s,m+\ell){\bf\Gamma}(-u+s,\ell+1-m)\, du
\nonumber \\ && \makebox[4em]{}
= \bigl( 2 \sh \pi(\alpha + \imath\beta) \bigr)^{-\ell-\frac12 - \imath s}
\bigl( 2 \sh \pi(\alpha - \imath \beta) \bigr)^{\ell+\frac12 - \imath s}.
\eea

The latter two formulas can be unified into one by redefining the range of values of the variables
$\ell$ and $m$. Namely, we set $m, \ell\in\Z+\mu,\, \mu=0, \tfrac12$. Then $2\ell$ is even for $\mu=0$ and
$2\ell$ is odd for $\mu=\tfrac12$. As a result, we can  write universally
\bea && \makebox[-2em]{}
\frac{1}{4\pi \bm{\Gamma}(2s, 2\ell)}
 \sum_{m \in \mathbb{Z}+\mu} e^{\pi\imath(m-\ell)} \int_{\mathbb{R}}
 e^{-2\pi \imath(  \alpha u +  \beta m) } \; {\bf\Gamma}(s\pm u,\ell \pm m)\, du
\nonumber \\ && \makebox[1em]{}
= \bigl( 2 \sh \pi(\alpha + \imath\beta) \bigr)^{-\ell - \imath s}
\bigl( 2 \sh \pi(\alpha - \imath \beta) \bigr)^{\ell - \imath s},\quad \ell\in\Z+\mu, \, \mu=0, \tfrac12.
\eea
At last, we shift the value of $\beta$, $\beta\to \beta+\tfrac12$, and come to our final formula
\bea \label{comp_red_pent} && \makebox[2em]{}
\frac{1}{4\pi}
 \sum_{m \in \mathbb{Z}+\mu}\int_{\mathbb{R}}
 e^{-2\pi \imath(  \alpha u +  \beta m) } \; {\bf\Gamma}(s\pm u,\ell \pm m)\, du
\nonumber \\ && \makebox[-1em]{}
=  \bm{\Gamma}(2s, 2\ell)\bigl( 2 \ch \pi(\alpha + \imath\beta) \bigr)^{-\ell - \imath s}
\bigl( 2 \ch \pi(\alpha - \imath \beta) \bigr)^{\ell - \imath s},\quad \ell\in\Z+\mu, \, \mu=0, \tfrac12.
\eea
We were assuming at the beginning that $\beta\in\mathbb{C}$. However, the sum of integrals
on the left-hand side of \eqref{qbin_comp} can converge only if $\beta\in\mathbb{R}$
and we can limit $\beta\in [0,1]$. For real $\beta$ the sum of integrals converges if
$\text{Im}(h)>-1$ (i.e. $\text{Im}(s)>-1/2$).

The derived beautiful formula coincides with the complex binomial theorem.
Indeed, shifting $\beta\to\beta+1/2$ in \eqref{qbin_comp} we can rewrite it as
\bea\label{qbin_comp_old} && \makebox[1em]{}
\frac{1}{4\pi \bm{\Gamma}(A|A')}
 \sum_{m \in \mathbb{Z}} \int_{\mathbb{R}-\imath 0}
  \frac{{\bf\Gamma}(s|s') {\bf\Gamma}(A-s,A'-s')}{Z^s \bar Z^{s'}}\, du= \frac{1}{(1+Z)^A (1+\bar Z)^{A'}},
\\ && \makebox[-1em]{}
Z=e^{2\pi(\alpha+\imath \beta)},\; \bar Z=e^{2\pi(\alpha-\imath \beta)},\quad \; A=\frac{r+\imath h}{2}, \; A-A'=r,\quad s=\frac{m+\imath u}{2},\; s-s'=m,
\nonumber\eea
which coincides with the formula given, e.g. in the Appendix of \cite{SaSpJGP}. Thus we come to an interesting
conclusion that the reduced pentagon identity \eqref{pentafkv5} {describes} a $q$-analogue
of the complex binomial theorem.

\section{The complex analogue of the Euler beta integral}

Consider three complex numbers $z, \alpha$ and $\alpha'$
such that $\alpha-\alpha' \in\mathbb{Z}$ and introduce the notation
$$
[z]^\alpha:= z^\alpha \bar z^{\alpha'}=(-1)^{\alpha-\alpha'}[-z]^\alpha,\quad
\int_{\mathbb{C}} d^2z:=\int_{\mathbb{R}^2}d(\text{Re}\, z)\, d(\text{Im}\, z),
$$
with $\bar z$ being the complex conjugate of $z$. Evidently, $[z]^\alpha$ is a single valued
function of $z$. Below we use this notation with a
reservation that in all formulas where $\alpha$ is an explicit integer, we set $\alpha'=\alpha$.
Then the complex analogue of the Euler beta integral
$$
\int_0^1 x^{a-1}(1-x)^{b-1}dx=\frac{\Gamma(a)\Gamma(b)}{\Gamma(a+b)},\quad \text{Re}(a), \text{Re(b)}>0,
$$
has the following form \cite{GGV}
\begin{equation}\label{compBF}
		\bm{\mathrm{B}}(a, b) = \frac{1}{\pi} \int_{\mathbb{C}}\, [t]^{a-1} [1-t]^{b-1}d^2t
= \frac{\bm{\Gamma}(a|a') \bm{\Gamma}(b|b')}{\bm{\Gamma}(a+b|a'+b')}.
\end{equation}
Equivalently, it can be rewritten as
\begin{equation}\label{CB}
\frac{1}{\pi} \int_{\mathbb{C}}[w-z_1]^{a-1} [z_2-w]^{b-1} d^2w
= \frac{{\bf\Gamma}(a,b,c)}{[z_1-z_2]^{c}},
\quad {\bf\Gamma}(a_1,\ldots,a_k):=\prod_{j=1}^k{\bf\Gamma}(a_j|a_j'),
\end{equation}
where $a+b+c= a'+b'+c'=1$. This integral converges for
\begin{equation}
	\Re(a + a') > 0, \qquad \Re(b + b') > 0, \qquad \Re(c+c') >0.
\end{equation}

Simple linear fractional transformations of the variables $w, z_1, z_2$ bring formula \eqref{CB}
to the following star-triangle relation form
\begin{eqnarray}
\int_{\mathbb{C}}[z_1-w]^{\alpha-1} [z_2-w]^{\beta-1} [z_3-w]^{\gamma-1}\frac{d^2w}{\pi}
=\frac{{\bf \Gamma}(\alpha,\beta,\gamma)} {[z_3-z_2]^{\alpha}[z_1-z_3]^{\beta}[z_2-z_1]^{\gamma}}.
 \label{STR}\end{eqnarray}

One can change the integration variable on the left-hand side of the exact evaluation formula \eqref{compBF} to
\begin{equation}
	t = \frac{s}{1 + s},
\end{equation}
which yields another standard form of the integral
\begin{align}\label{compBF2}
	\bm{\mathrm{B}}(a,b) = \frac{1}{\pi} \int_{\mathbb{C}}\frac{[s]^{a-1}}{[1 + s]^{a+b}}\, d^2s.
\end{align}
For  $z_1 = 0,\, z_2 = -1$ the right-hand side of relation  \eqref{CB} is symmetric in $a, b, c$.
Therefore one can permute these parameters on the left-hand side. The replacement $b\to 1-a-b$ yields the
expression \eqref{compBF2} up to the sign factor $(-1)^{a-a'+b-b'}$.

We change the integration variable in  \eqref{compBF2}
\begin{align}\label{s_change}
s = e^{\gamma},\qquad	\gamma = \alpha + \imath \beta \in \mathcal{C} \colon \quad \alpha \in \mathbb{R}, \quad \beta \in [0, 2\pi)
\end{align}
and denote
\begin{align}
	a =  z+ \frac{g}{2}, \quad a' = z' + \frac{g'}{2}, \quad b = -  z+ \frac{g}{2}, \quad b' =  -z' + \frac{g'}{2},
\end{align}
so that $a+b=g,\, a'+b'=g'$. Define $\mathcal{L}_\mu$ as a set of variables
$$
x=\frac{k + \imath u}{2} \in \mathcal{L}_\mu \colon \quad k \in \mathbb{Z}+\mu,\;\mu=0,\, \tfrac12,
\quad u \in \mathbb{C}.
$$
If we denote $x'=\tfrac{-k+\imath u}{2}$, then $x-x'=k\in\Z +\mu$, i.e. for $\mu=\tfrac12$ this
difference is not an integer. Therefore this notation should be used with a caution that
for $\mu=\tfrac12$ only for the sums (or differences) $x+y$ of pair elements $x, y\in\mathcal{L}_\mu$
one has $(x+y)-(x+y)'\in\Z$. The notations used below should be understood with this reservation.

The convergence conditions of the beta integral take now the form
$$
0<\text{Re}(g+g')<2,\quad  -\tfrac12 \text{Re}(g+g')< \text{Re}(z+z')< \tfrac12 \text{Re}(g+g').
$$
Since we should have $a-a', b-b'\in\Z$, in general $z,\, \tfrac12 g \in \mathcal{L}_\mu$.
In particular, in the representation $z=\tfrac12(k+\imath u)$ the discrete variable $k$ may be
either integer ($\mu=0$) or half-integer ($\mu=\tfrac12$).
As a result of the change \eqref{s_change}, the integration measure $\int_{\mathbb C}d^2s$ changes to $\int_{\mathcal{C}}e^{2\alpha}d^2\gamma$ and we obtain
\begin{equation}
\bm{\mathrm{B}}(a,b)=\frac{1}{\pi}\int_{-\infty}^\infty d\alpha\int_0^{2\pi}d\beta
\frac{e^{\imath(\alpha u+\beta k)}}{( 2 \ch \frac{\alpha+\imath\beta}{2})^{g}( 2 \ch \frac{\alpha-\imath\beta}{2})^{g'}}
=\frac{\bm{\Gamma}\bigl(\pm  z+ \tfrac{g}{2} \big| \pm z' + \tfrac{g'}{2} \bigr)}{\bm{\Gamma}(g|g')}.
\label{comp_beta_mod}\end{equation}

Suppose now that $u$ is a real variable, $u\in\R$. Then $z'=-\bar z$ and this formula
represents the two-dimensional Fourier transformation of the powers of elementary functions
(with the reservation that $k=n+\mu,\, n\in\Z,\, \mu=0,\,\tfrac12$).
\begin{equation}
\bm{\mathrm{B}}(a,b)=\frac{1}{\pi}\int_{\mathcal C}
\frac{e^{\gamma z-\bar\gamma\bar z}}{[ 2 \ch \frac{\gamma}{2}]^{g}} \,d^2\gamma
=\frac{\bm{\Gamma}\bigl(\pm  z+ \tfrac{g}{2} \big| \mp \bar z + \tfrac{g'}{2} \bigr)}{\bm{\Gamma}(g|g')}.
\label{comp_beta_mod2}\end{equation}
Performing the inverse Fourier transformation we obtain the equality
\begin{equation}
\frac{1}{4\pi{\bm{\Gamma}(g|g')}}\int_{\mathcal{L}_\mu }
e^{-\gamma z+\bar\gamma\bar z}\,
{\bm{\Gamma}\bigl(\pm  z+ \tfrac{g}{2} \big| \mp \bar z + \tfrac{g'}{2} \bigr)}\,d^2z
=\Big[2 \ch \frac{\gamma}{2}\Big]^{-g},
\label{comp_beta_mod3}\end{equation}
where $\int_{\mathcal{L}_\mu}d^2z:=\sum_{k\in\Z+\mu}\int_{-\infty}^{\infty}du.$
After replacing $\alpha, \,\beta$ by $2\pi\alpha,\, 2\pi\beta$,
we obtain exactly the complex binomial theorem \eqref{comp_red_pent}, \eqref{qbin_comp_old}.
So, the complex beta integral is basically nothing else than the analytically continued
inverse Fourier transform of the complex binomial theorem.

Note that the relation  \eqref{pentafkv5} is evidently self-dual with respect to the Fourier transform.
Therefore there should exist some direct degeneration limit from  \eqref{pentafkv5} to \eqref{comp_beta_mod2}.
Indeed, consider the hyperbolic beta integral \eqref{pentafkv5} in a different parametrization
corresponding to the inverse Fourier transformation
\begin{align}\label{pent_inv}
	\int_{-\imath\infty}^{\imath\infty} e^{\frac{2\pi \imath \lambda z}{\omega_1 \omega_2}}  \gamma^{(2)}\biggl( \pm \lambda + \frac{g}{2} ; \omega \biggr) \, \frac{d\lambda}{\imath \sqrt{\omega_1 \omega_2}} = \frac{1}{\gamma^{(2)}(\bar{g}; \omega)}  \gamma^{(2)}\biggl( \pm z + \frac{\bar{g}}{2} ; \omega \biggr),
\end{align}
where $\bar{g} = \omega_1 + \omega_2 - g$. Fix, as before, $\omega_1 = \imath + \delta,
\, \omega_2 = - \imath + \delta,$ parametrize
\begin{align}
	& z = \imath \sqrt{\omega_1 \omega_2} (m + u \delta), && m \in \mathbb{Z}, \qquad u \in \mathbb{C}, \\[6pt]
	& \bar{g} = \imath \sqrt{\omega_1 \omega_2} (2\ell + 2s \delta), && \ell \in \mathbb{Z}, \quad s \in
\mathbb{C},
\end{align}
and take the limit $\delta \to 0^+$. For simplicity, we restricted description to the case $\mu=0$ when $m\in\Z$. Then the asymptotics of the right-hand side in \eqref{pent_inv} is given by the expression
\begin{align}\label{rhs-lim-beta}
	 \frac{1}{\gamma^{(2)}(\bar{g}; \omega)}  \gamma^{(2)}\biggl( \pm z + \frac{\bar{g}}{2} ; \omega \biggr) \underset{\delta \to 0^+}{=} e^{\pi \imath(m - \ell)} \frac{\bm{\Gamma}(s \pm u, \ell \pm m)}{4\pi \delta \, \bm{\Gamma}(2s, 2\ell)},
\end{align}
which is equal up to the diverging factor $1/\delta$ to the integrand in formula \eqref{even_l}
with the exclusion of the Fourier exponential
$e^{-2\pi \imath(  \alpha u +  \beta m) }$. This fact already indicates that the
left-hand side expression in \eqref{pent_inv} should also diverge in an analogous way and be
proportional to the left-hand side of \eqref{comp_beta_mod}. We were able to find a generalization
of relation \eqref{intdel} that  leads to the limit
\begin{align}
\stackreb{\lim}{\delta\to 0^+}
& \delta\int_{-\imath\infty}^{\imath\infty} e^{\frac{2\pi \imath \lambda z}{\omega_1 \omega_2}}  \gamma^{(2)}\biggl( \pm \lambda + \frac{g}{2} ; \omega \biggr) \, \frac{d\lambda}{\imath \sqrt{\omega_1 \omega_2}}
\\ \nonumber & \; =
\int_{-\frac{1}{2}}^{\frac{1}{2}} d\beta \int_{\mathbb{R}} e^{-2\pi \imath(\alpha u + \beta m)} \bigl( 2\sh \pi(\alpha + \imath \beta ) \bigr)^{-\ell - \imath s} \, \bigl(2 \sh \pi(\alpha - \imath \beta) \bigr)^{\ell - \imath s} d\alpha.
 \end{align}
Combining it with the right-hand side asymptotic expression~\eqref{rhs-lim-beta}, we arrive at the desired identity
\begin{eqnarray*} &&
	\int_{-\frac{1}{2}}^{\frac{1}{2}} d\beta \int_{-\infty}^\infty e^{-2\pi \imath(\alpha u + \beta m)} \bigl( 2\sh \pi(\alpha + \imath \beta ) \bigr)^{-\ell - \imath s} \, \bigl(2 \sh \pi(\alpha - \imath \beta) \bigr)^{\ell - \imath s} d\alpha
\\  && \makebox[3em]{}
	= e^{\pi \imath(m - \ell)} \frac{\bm{\Gamma}(s \pm u, \ell \pm m)}{4\pi \, \bm{\Gamma}(2s, 2\ell)},
\end{eqnarray*}
which coincides with the complex beta integral evaluation formula \eqref{comp_beta_mod} after shifting
$\beta\to\beta-1/2$ and rescaling of $\alpha$ and $\beta$. The details of this consideration, which looks like
a combined application of the method of \cite{Sarkissian:2020ipg} together with the quasiclassical limit
considered in \cite{DM2019}, will be presented later on.

As shown in \cite{DMV2017}  with the help of the complex
binomial theorem \eqref{qbin_comp_old}, the equality \eqref{limtrafoI}
and thus the complex pentagon identity \eqref{pent_comp} {represent}
the Mellin-Barnes form of the star-triangle relation \eqref{STR}
(i.e., effectively, of the complex analogue of the Euler beta integral).
Moreover,  relation  \eqref{STR} is easily proved with the help of  \eqref{qbin_comp_old}
and complex beta integral, see Appendix of \cite{SaSpJGP}.
This means that the Faddeev-Volkov star-triangle relation \eqref{hyper2},
being a $q$-analogue of \eqref{limtrafoI} (or of \eqref{STR}), also should be
directly related to \eqref{pentafkv5}. Indeed, this relation was noticed in \cite{Kashaev,Volkov}
and we give here a brief outline of this fact in our notation.

Let us rewrite the left-hand side of equality \eqref{hyp4-5} in the form
\begin{equation}\label{45}
\int_{-\textup{i}\infty}^{\textup{i}\infty}
e^{\frac{\pi \textup{i} z}{\omega_1\omega_2}(f+h+g)}
\gamma^{(2)}(f+z,h-z,\tfrac{g}{2}\pm z;\mathbf{\omega})
{dz\over \textup{i}\sqrt{\omega_1\omega_2}},
\end{equation}
where we denoted $f=f_1, h=h_1, f_2=h_2=\tfrac{g}{2}$ (this can be reached after appropriate shift of the
integration variable $z$ and removal of the arising additional exponent).
Now we  replace in formula \eqref{pentafkv5} the integration variable $z$ by $y$ and $\lambda$ by $z$.
Then we take the product $\gamma^{(2)}(\tfrac{g}{2}\pm z;\mathbf{\omega})$ from the right-hand
side and substitute it into the above expression. This yields
$$
\gamma^{(2)}(g;\mathbf{\omega})
\int_{(\textup{i}\R)^2}
e^{2\pi\textup{i}z(y+\frac12 (f+h+g)\over \omega_1\omega_2}
\gamma^{(2)}(f+z,h-z,\pm y+{\bar{g}\over 2};\mathbf{\omega})
{dy\over \textup{i}\sqrt{\omega_1\omega_2}}
{dz\over \textup{i}\sqrt{\omega_1\omega_2}}
$$
Now we integrate over the variable $z$ using the same identity  \eqref{pentafkv5} after shifting $z\to z+\tfrac{h-f}{2}$. This leads to the expression
$$
\gamma^{(2)}(g,f+h;\mathbf{\omega})\int_{\textup{i}\R}
e^{\pi\textup{i}(h-f)(y+\frac12(f+h+g))\over \omega_1\omega_2}
\gamma^{(2)}(y+{\bar{g}\over 2},-y+\frac{\bar g}{2}-f-h;\mathbf{\omega})
{dy\over \textup{i}\sqrt{\omega_1\omega_2}},
$$
where we used the identity $\gamma^{(2)}(y+\tfrac{\omega_1+\omega_2+g}{2},-y+\tfrac{\bar{g}}{2};\mathbf{\omega})=1$.
Finally, we shift the integration variable $y\to y-\tfrac{f+h}{2}$ and using once again relation
\eqref{pentafkv5} we obtain the final result
$$
e^{\frac{\pi\textup{i}(h-f)g}{2\omega_1\omega_2}}\gamma^{(2)}(g,f+h,\bar g-f-h,\tfrac{g}{2}+f,
\tfrac{g}{2}+h;\mathbf{\omega}),
$$
which coincides with the right-hand side of equality  \eqref{hyp4-5}.

As to the general relation \eqref{hyper2}, let us denote as $I$ the corresponding integral. Apply
to the product $\gamma^{(2)}(f_3+z,h_3-z;\mathbf{\omega})$ in the integrand the above trick. This yields
\begin{eqnarray}\label{fl} &&
I=\gamma^{(2)}(f_3+h_3;\mathbf{\omega}) \int_{\imath\R}e^{\frac{\pi\imath x}{\omega_1\omega_2} (f_3-h_3)}
 \gamma^{(2)}(\tfrac{\omega_0 -f_3-h_3}{2}\pm x;\mathbf{\omega})
 F(\underline{f}; \underline{g};x){dx\over \textup{i}\sqrt{\omega_1\omega_2}},
\end{eqnarray}
where $\omega_0=\omega_1+\omega_2$ and
\be\label{newf}
F(\underline{f}; \underline{g};x)=\int_{\imath\R}e^{{2\pi\textup{i}xz\over \omega_1\omega_2}}\prod_{a=1,2}\gamma^{(2)}(f_a+z,h_a-z;\mathbf{\omega}){dz\over \textup{i}\sqrt{\omega_1\omega_2}}.
\ee
Function $F$ looks similar to \eqref{45}, but now the Fourier exponential has a free parameter.

At the next step the pentagon identity \eqref{pentafkv5} generates symmetry transformation of function
\eqref{newf}. To derive that it is necessary to apply the same trick as above to the products of gamma functions
in the integrand of \eqref{newf} separately for parameters with indices $a=1$ and $a=2$ and integrate over $z$.
This yields a $\delta$-function which is eliminated by one of the integrations. As a result,
one obtains a dual representation of \eqref{newf}
\bea\nonumber
&&F(f_1,f_2;h_1,h_2;x)=e^{-\frac{\pi \textup{i}x}{\omega_1\omega_2}\beta_2} \gamma^{(2)}(\alpha_1,\alpha_2;\mathbf{\omega})
\\  && \makebox[2em]{}
\times F\left(\tfrac{\omega_0 -\alpha_1}{2},x+\tfrac{\omega_0 -\alpha_2}{2};
\tfrac{\omega_0 -\alpha_1}{2},-x+\tfrac{\omega_0 -\alpha_2}{2};\tfrac{\beta_1-\beta_2}{2}\right),
\label{trafo}\eea
where we denoted
$$
\alpha_a=f_a+h_a, \quad \sum_{a=1}^3\alpha_a=\omega_0,\quad \beta_1=f_1-h_1,\quad \beta_2=f_2-h_2.
$$

The derived symmetry transformation can be iterated after permuting arguments of the $F$-function.
After permutation of the first two arguments and applying transformation \eqref{trafo}
one obtains a relation similar to \eqref{trafo} with four $\gamma^{(2)}$-functions in front of the
$F$-function on the right-hand side. Repeating this procedure one more time (i.e. permute the two
first arguments of the resulting $F$-function and apply again \eqref{trafo}), one comes to the identity
\begin{eqnarray}\nonumber &&
F(f_1,f_2;h_1,h_2;x)=e^{\frac{\pi \textup{i}x}{\omega_1\omega_2} (h_1-f_2)}
\gamma^{(2)}(f_1+h_1,f_2+h_2,f_1+h_2,f_2+h_1;\mathbf{\omega})
\\ \label{find1} && \makebox[2em]{}
\times\gamma^{(2)}(\tfrac{\omega_0+f_3+h_3}{2}\pm x;\mathbf{\omega})
\\ \nonumber && \makebox[2em]{}
\times F\big(\tfrac{\omega_0-f_1-h_2}{2},\tfrac{\omega_0+f_1-h_2}{2}-f_2;
\tfrac{\omega_0-f_1-h_2}{2},\tfrac{\omega_0-f_1+h_2}{2}-h_1;x\big).
\end{eqnarray}
As a result of these manipulations, most importantly, the variable $x$ disappeared
from the first four arguments of the $F$-function, i.e. it enters now only in
the exponential in the definition \eqref{newf}.

At the final step, one should substitute the right-hand side expression of the equality
\eqref{find1} into \eqref{fl}. After that four $\gamma^{(2)}$-functions
get cancelled and the integral over $x$ brings in the
delta-function $\delta(z+{1\over 2}(f_3-h_3-f_2+h_1)$. Eliminating this delta-function by integrating over
$z$, one arrives to the right-hand side expression of the general identity \eqref{hyper2}.
As expected, only the pentagonal relation \eqref{pentafkv5} (a hyperbolic analogue of the complex
binomial theorem) was used in the proof. However, this proof is much longer than the
proof of the complex star-triangle relation with the help of the complex binomial theorem \cite{SaSpJGP}.

\section{Degeneration of the Fourier transformation of the hyperbolic gamma function}

Consider the simplest hyperbolic hypergeometric integral --- the Fourier transformation of the modular dilogarithm \eqref{namer3}
\begin{eqnarray*} &&
\int_{\textup{i}\mathbb{R}+0}
e^{\frac{2\pi \textup{i}\lambda z}{\omega_1\omega_2} + \frac{\pi \textup{i}}{2}(B_{2,2}(z;\mathbf{\omega})
- B_{2,2}(0;\mathbf{\omega}))}
\gamma^{(2)}(z;\mathbf{\omega})\frac{dz}{\textup{i}\sqrt{\omega_1\omega_2}}
=e^{-{\pi \textup{i}\over 2} B_{2,2}(\lambda;\mathbf{\omega})}
\gamma^{(2)}(\lambda;\mathbf{\omega}).
\end{eqnarray*}
Let us argue that its degeneration to the complex hypergeometric level
in the limit $b\to\imath$ \eqref{om1om22} is given by the formula
\begin{align}\label{int2}
\frac{1}{4\pi}\sum_{k \in \mathbb{Z}} \int_{\mathbb{R}-\imath 0}  e^{- \imath \alpha u - \imath \beta k} \; \bm{\Gamma}(u, k)\, du  = \exp( e^{\alpha - \imath \beta} - e^{\alpha + \imath \beta}),
\end{align}
which is the two-dimensional Fourier transform of the complex gamma function. {For $\alpha, \beta \in \mathbb{R}$ we have conditional convergence of the integral in $u$ and absolute convergence of the series in $k$, see \cite[Section 1.4]{Neretin2019}. }

Following the limiting formula for the hyperbolic gamma function \eqref{gam2lim2},
we parametrize the integration variable $z$ in \eqref{namer3} correspondingly
$z = \imath \sqrt{\omega_1 \omega_2} (k + u \delta)$
and take the limit $\delta\to 0^+$. As before, the poles of $\gamma^{(2)}(z;\mathbf{\omega})$-function
approach the integration contour and the integrand blows up. However, this time there is no pinching
of the integration contour, since the poles are located only from the left side of the contour.
The divergency is removed by rescaling the integration
variable, $z\propto \delta u$, and we obtain a finite limit in the form of an infinite sum of
integrals of complex gamma functions over $u$.

In order to cancel the diverging factor $(4 \pi \delta)^{\imath u}$ in the integrand of \eqref{namer3},
we take the parameter $\lambda$ in the form
$$
	\lambda = \tilde\lambda + \imath \sqrt{\omega_1 \omega_2} \, \frac{\ln (4\pi \delta_N)}{2\pi \delta_N},
$$
where $\tilde \lambda$ will be specified below and $\delta \equiv \delta_N$ is the function of $N \in \mathbb{Z}_{>0}$ defined by the equation
\begin{align}
	\frac{\ln (4\pi \delta_N)}{4 \pi \delta_N} = - N.
\end{align}
It can be expressed in terms of the Lambert $W(x)$ function
\begin{align}
	\delta_N = \frac{e^{-W(N)}}{4\pi}, \qquad  W(x) \, e^{W(x)} = x,
\end{align}
which is uniquely defined for $x \geq 0$. Moreover, it has the asymptotics
\begin{align}
	W(x) \underset{x \to \infty}{=} \ln x - \ln (\ln x),
\end{align}
so that
\begin{align}\label{d-as}
	\delta_N \underset{N \to \infty}{=} \frac{\ln N}{4\pi N} \to 0^+.
\end{align}
For such a choice of $z$ and $\lambda$ the exponent in the integrand takes the form
$$
	e^{\frac{2\pi \imath \lambda z}{\omega_1 \omega_2} } = e^{\frac{2\pi \imath z\tilde\lambda}{\omega_1 \omega_2}  + 4\pi \imath N k -\imath u\log (4\pi \delta_N)}
 = e^{\frac{2\pi \imath z\tilde\lambda}{\omega_1 \omega_2}} \; (4\pi \delta_N)^{- \imath u}  .
$$

It is left to specify $\tilde\lambda$ in such a way that in the limit $\delta_N\to 0^+$ both, the value of
$\delta_N\tilde \lambda$ and the exponential $\exp(-\frac{2\pi k \tilde\lambda}{\sqrt{\omega_1 \omega_2}})$,
remain finite. {With such an input, we consider
\begin{align}
\tilde\lambda = \delta_N+\imath \sqrt{\omega_1 \omega_2} \biggl( M_N + \frac{\beta}{2\pi} \biggr), \qquad M_N = \biggl\lfloor \frac{2 \alpha N}{\ln N} \biggr\rfloor \in \mathbb{Z},
\end{align}
where the first term was added in order to have Re$(\lambda)>0$. Here
}
$\alpha, \beta \in \mathbb{R}$ and $\lfloor x \rfloor$ denotes the integer part of $x$. In this case
\begin{align}
	M_N \, \delta_N \underset{N \to \infty}{=} \frac{\alpha}{2\pi}.
\end{align}
Indeed,
\begin{align}
	M_N \, \delta_N = \frac{\alpha}{2\pi} \, \frac{4\pi N}{\ln N} \, \delta_N + \biggl(  \biggl\lfloor \frac{2 \alpha N}{\ln N} \biggr\rfloor - \frac{2\alpha N}{\ln N} \biggr) \delta_N.
\end{align}
The first term tends to $\alpha/ 2\pi$ due to the asymptotics of $\delta_N$ \eqref{d-as}.
The second term tends to zero, since the expression in the brackets is bounded and
$\delta_N \to 0^+$. As a result, the parametrization
\begin{align}
	z = \imath \sqrt{\omega_1 \omega_2}(k + u \delta_N), \quad \lambda = \imath \sqrt{\omega_1 \omega_2} \biggl( \biggl\lfloor \frac{2 \alpha N}{\ln N} \biggr\rfloor + \frac{\beta}{2\pi} + \frac{\ln(4\pi \delta_N)}{2\pi \delta_N} \biggr)
\end{align}
leads to the limiting relation
\begin{align}
	e^{\frac{2\pi \imath z\lambda}{\omega_1 \omega_2} }
 \underset{N \to \infty}{=} e^{- \imath \alpha u - \imath \beta k}.
\end{align}

The exponential of the Bernoulli polynomials in \eqref{namer3}  $e^{\frac{\pi \textup{i}}{2}(B_{2,2}(z;\mathbf{\omega})- B_{2,2}(0;\mathbf{\omega}))}$ in the limit $\delta_N\to 0$
simplifies to $e^{-\frac{\pi\imath}{2}k^2}$ which cancels a similar factor arising in
the asymptotics \eqref{gam2lim2}. As a result, the left-hand side of \eqref{namer3} has the limit
\begin{eqnarray} \nonumber &&
\int_{\textup{i}\mathbb{R}+0}
e^{\frac{2\pi \textup{i}\lambda z}{\omega_1\omega_2} + \frac{\pi \textup{i}}{2}(B_{2,2}(z;\mathbf{\omega})
- B_{2,2}(0;\mathbf{\omega}))}
\gamma^{(2)}(z;\mathbf{\omega})\frac{dz}{\textup{i}\sqrt{\omega_1\omega_2}}
\\ && \makebox[2em]{} \label{lhs-lim}
	\underset{N \to \infty}{=} \;
\frac{1}{4\pi}\sum_{k \in \mathbb{Z}} \int_{\mathbb{R}-\imath 0}  e^{- \imath \alpha u - \imath \beta k} \; \bm{\Gamma}(u, k)\, du.
\label{limit}\end{eqnarray}

On the right-hand side of identity \eqref{namer3} we have
$$
	e^{- \frac{\pi \imath}{2} B_{2,2}(\lambda; \omega) } \; \gamma^{(2)}(\lambda; \omega)
= \frac{e^{- \frac{\pi \imath}{2} B_{2,2}(\lambda; \omega) } }{\gamma^{(2)}(\omega_1 + \omega_2 - \lambda; \omega)}
= \frac{\bigl( e^{2\pi \imath\frac{\omega_1 - \lambda}{\omega_2} }; q
\bigr)_{\! \infty}}{\bigl( e^{-2\pi \imath\frac{\lambda}{\omega_1} }; \tilde q
\bigr)_{\! \infty} }.
$$
This is a ratio of two $q$-exponential functions \cite{GR}
\begin{align}\nonumber
	{}_1 \varphi_0(0;q, z) = \sum_{k = 0}^\infty \frac{z^n}{(q;q)_n} = \frac{1}{(z; q)_{\infty}}, \qquad |z| < 1,
\end{align}
with the well known limiting relation
\begin{align}\label{q-exp_lim}
{}_1 \varphi_0(0;q,(1-q)z) = \frac{1}{((1-q)z; q)_{\infty}} \underset{q\to 1^-}{=} e^z.
\end{align}
With the parametrization of $\omega_1, \omega_2$ and $\lambda$ taken above, in the limit $N\to\infty$
we have
\begin{align}
	& e^{2\pi \imath\frac{\omega_1 -  \lambda}{\omega_2}} \underset{N \to \infty}{=} 4\pi \delta_N \, e^{\alpha + \imath \beta+O(\delta_N\ln\delta_N)}, && q=e^{2\pi \imath \frac{\omega_1}{\omega_2} } \underset{N \to \infty}{=} e^{- 4\pi \delta_N+O(\delta_N^2)} \\[6pt]
	& e^{-2\pi \imath\frac{\lambda}{\omega_1} } \underset{N \to \infty}{=} 4\pi \delta_N \, e^{\alpha - \imath \beta+O(\delta_N\ln\delta_N)}, && \tilde q=e^{- 2\pi \imath \frac{\omega_2}{\omega_1} } \underset{N \to \infty}{=} e^{- 4\pi \delta_N+O(\delta_N^2)}.
\end{align}
Using the limiting relation \eqref{q-exp_lim} when both $1-q$ and $1-\tilde q$ are approximated as
$4\pi \delta_N+O(\delta_N^2)$, we find
$$
\Bigl( e^{2\pi \imath\frac{\omega_1 - \lambda}{\omega_2} }; q\Bigr)_{\! \infty} \; \underset{N \to \infty}{=} \; \exp \bigl( - e^{\alpha + \imath \beta} \bigr),
\quad \Bigl( e^{-2\pi \imath\frac{\lambda}{\omega_1} }; \tilde q \Bigr)_{\! \infty} \; \underset{N \to \infty}{=} \; \exp \bigl( - e^{\alpha - \imath \beta} \bigr).
$$
Finally, we obtain for the right-hand side of \eqref{namer3}
\begin{align}\label{rhs-lim}
	e^{- \frac{\pi \imath}{2} B_{2,2}( \lambda; \omega) } \; \gamma^{(2)}( \lambda; \omega) \underset{N \to \infty}{=}  \exp( e^{\alpha - \imath \beta} - e^{\alpha + \imath \beta}).
\end{align}
Equating the limits of both sides in equality \eqref{namer3}, we come to the stated relation \eqref{int2}.

Let us take the inverse Fourier transformation of the derived identity. It yields
\begin{eqnarray} \nonumber && % \label{2comp_exp} &&
	\bm{\Gamma}(u, k) = \frac{1}{\pi} \int_{\mathbb{R}} d\alpha \int_{0}^{2\pi} \! d\beta
\; e^{\imath \alpha u + \imath \beta k} \, \exp(e^{\alpha - \imath \beta} - e^{\alpha + \imath \beta})
\\ && \makebox[3em]{}
	= \frac{1}{\pi} \int_{\mathbb{C}} z^{a - 1 }
\, \bar{z}^{a' - 1 } \; e^{\bar{z} - z}\, d^2 z=\bm{\Gamma}(a|a'), \quad a=   \frac{k + \imath u}{2}.
 \label{GGR}\end{eqnarray}
In order to obtain the second representation we changed the integration variable $z = e^{\alpha + \imath \beta}$
in the Fourier transformation of the double exponential function
and passed to the integration over complex plane. Formula \eqref{GGR} coincides with the one
used in \cite{GGR} for defining the complex gamma function.
So, the formula we obtained \eqref{int2} is {a mere} inverse Fourier transform of the latter integral.

Representation of the complex gamma function \eqref{GGR} provides a simple way to evaluate the complex
analogue of the Euler beta integral. Indeed, consider the double complex integral
$$
\bm{\Gamma}(a|a')\bm{\Gamma}(b|b')= \frac{1}{\pi^2} \int_{\mathbb{C}^2} z^{a - 1 }
\, \bar{z}^{a' - 1 } \; e^{\bar{z} - z} w^{b - 1 }
\, \bar{w}^{b' - 1 } \; e^{\bar{w} - w}\, d^2z\, d^2 w.
$$
Change the integration variables
$$
u=z+w,\quad z=uv,
$$
which gives $w=u(1-v)$ and $d^2z\,d^2w=|u|^2d^2u\,d^2v$. As a result, we come to the integral
$$
\frac{1}{\pi^2} \int_{\mathbb{C}} [u]^{a+b - 1 } e^{\bar{u}-u}\, d^2 u
\int_{\mathbb{C}} [v]^{a-1}[1-v]^{b-1}\, d^2v,
$$
which immediately leads to formula \eqref{compBF}. In this picture the main complexity goes to
the evaluation of integral \eqref{GGR} as the ratio of two Euler gamma functions which is a
nontrivial task \cite{GGR}.

We also remark that the limit \eqref{rhs-lim} is a complex analogue of the formula
(A.33) from \cite{RuijSIGMA}, which in our notation reads
\begin{align}
	e^{- \frac{\pi \imath}{2} B_{2,2}( \lambda; \omega) } \; \gamma^{(2)}( \lambda; \omega) \underset{\omega_1 \to 0^+}{=} \exp (- e^{\alpha}),
\end{align}
where
\begin{align}
	\lambda = \frac{\imath \omega_2}{2\pi} \biggl( \alpha + \ln \frac{2\pi \omega_1}{\omega_2} \biggr) + \frac{\omega_1 + \omega_2}{4}.
\end{align}
Using this limit together with the reduction to the ordinary gamma function \cite{Ruij}
\begin{align}
	\gamma^{(2)}( x \omega_1 ; \omega) \underset{\omega_1 \to 0^+}{=} \frac{1}{\sqrt{2\pi}} \biggl( \frac{2\pi \omega_1}{\omega_2} \biggr)^{x - \frac{1}{2}} \, \Gamma(x),
\end{align}
one can obtain from \eqref{namer3}
\begin{align}
\frac{1}{2\pi}\int_{\mathbb{R}-\imath 0}  e^{- \imath \alpha u} \, \Gamma(\imath u)\, du= \exp(-e^{\alpha}).
\end{align}
This identity is a real analogue of the formula \eqref{int2}.

\section{Conclusion}

Thus, after {a substantial} generalization of the method of degenerating hyperbolic hypergeometric integrals
rigorously justified in  \cite{Sarkissian:2020ipg}, we have derived the degeneration limits of the
simplest hyperbolic identities involving the exponential functions to the complex hypergeometric level.
This allowed us to clarify deep
relations between the complex binomial theorem, complex analogue of the Euler beta integral evaluation
and various pentagon identities. As a next step, we plan to apply the derived results to mathematical
physics problems. For instance, they are needed for consideration of the complex versions of integrable
many-body systems {found in} \cite{SaSpJPA}. In particular, they provide {necessary tools for}
an appropriate degeneration limit
for the wave functions of two-body Ruijsenaars hyperbolic system \cite{RuijSIGMA}, which was
recently reexamined in detail in \cite{BDKK2} in view of the solutions of eigenvalue problems
for general $N$-body system Hamiltonians \cite{BDKK1,HR}.

\smallskip

{\bf Acknowledgments.}
This study has been partially supported by the Russian Science Foundation (grant 24-21-00466).
Main results of this paper were obtained during N. Belousov's work at the Laboratory of Mirror Symmetry
and Automorphic Forms of NRU HSE.
\vspace*{0.5em}

%\newpage

\end{document}